\documentclass{amsart}
\usepackage{amssymb}
\usepackage{amsmath}
\usepackage{color}
\usepackage{graphicx}

\newtheorem{theorem}{Theorem}[section]
\newtheorem{lemma}[theorem]{Lemma}

\theoremstyle{definition}
\newtheorem{definition}[theorem]{Definition}

\theoremstyle{remark}
\newtheorem{remark}[theorem]{Remark}

\DeclareMathOperator*{\esssup}{ess\,sup}
\newcommand{\ds}{{\ell}^{p_n}(\mathbb{Z})}

\numberwithin{equation}{section}

\begin{document}

\title[Boundedness of discrete Hilbert transform and Discrete Mikhlin Multiplier on discrete variable Lebesgue space ]{discrete Hilbert transform and Discrete Mikhlin Multiplier on discrete variable Lebesgue space }

\author{Arash Ghorbanalizadeh}
\address{Department of Mathematics, Institute for Advanced Studies in Basic Sciences (IASBS), Zanjan 45137-66731, Iran}
\curraddr{}
\email{ghorbanalizadeh@iasbs.ac.ir}

\author{Reza Roohi Seraji}
\address{
Department of Mathematics, Institute for Advanced Studies in Basic Sciences (IASBS), Zanjan 45137-66731, Iran}
\curraddr{}
\email{rroohi@iasbs.ac.ir}


\date{}

\dedicatory{}

\commby{}


\begin{abstract}
In this paper, by using continuous Hilbert transform and maximal operator boundedness property in the variable Lebesgue space $ L^{p(\cdot)}(\mathbb{R}) $ we show that the discrete Hilbert transform is bounded in the variable discrete Lebesgue space $ \ell^{p_n}(\mathbb{Z}) $. We show that the discrete Mikhlin multiplier $ \mathcal{T}_m  $
is a bounded operator on $ \ell^{p_n}(\mathbb{Z})  $ when $ 1<\underline{p}_n<\bar{p}_n<\infty $.
\end{abstract}

\maketitle

\section{Introduction}
\subsection{Variable Lebesgue Spaces}
Variable Lebesgue spaces are a generalization of the classical Lebesgue spaces. These spaces were originally introduced by Orlicz \cite{WO} and are denoted by $L^{p(\cdot)}$.

Let us  revisit the notion of variable exponent Lebesgue spaces. We denote by $\mathcal{P}_0$ the set of measurable functions $p(\cdot): \mathbb{R}^n \rightarrow [c,\infty]$, where $c>0$, and by $\mathcal{P}$ the subset of $\mathcal{P}_0$ such that the range of its elements is contained in $[1,\infty]$.

Let $\Omega_{\infty}=\{x\in \mathbb{R}^n: p(x)=\infty\}$ and $\Omega_0=\mathbb{R}^n\setminus \Omega_{\infty}$ be the subsets of $\mathbb{R}^n$ where the exponent is infinite and finite, respectively. The variable exponent Lebesgue space $L^{p(\cdot)}(\mathbb{R}^n)$ consists of all measurable functions $f$ such that there exists a positive constant $\lambda$ for which the modular

\[\varrho_{L^{p(\cdot)}}(\lambda^{-1}f):=\int_{\Omega_0}\left(\frac{|f(x)|}{\lambda}\right)^{p(x)}
dx
+\left\|\lambda^{-1}f\right\|_{L^{\infty}(\Omega_{\infty})}
\]
is finite. Given a variable exponent  $p(\cdot)\in \mathcal{P}_0$   and a function $f \in L^{p(\cdot)}(\mathbb{R}^n)$, we introduce the Luxemburg norm associated with the variable exponent Lebesgue space as
\[
\|f\|_{L^{p(\cdot)}}=\inf\{\lambda>0:  \varrho_{L^{p(\cdot)}}(\lambda^{-1}f)\leq 1\}.
\]
There are two important continuity condition on $ p(\cdot) $ as follows \cite{Cruz},
\begin{itemize}
\item
We say that $ p(\cdot) $ is locally log-H\"{o}lder, and we write $ p(\cdot)\in LH_0(\mathbb{R}^n) $, if
\begin{align*}
\exists C_0\quad:\quad |x-y|<\frac{1}{2}\longrightarrow~|p(x)-p(y)|\leq\frac{C_0}{-\log (|x-y|)}
\end{align*}

\item
We say that $ p(\cdot) $ is log-H\"{o}lder continuous at infinity, and we write $ p(\cdot)\in LH_\infty(\mathbb{R}^n) $, if
\begin{align*}
\exists C_0,~\exists p_\infty\quad:\quad \forall x\in\mathbb{R}^n\longrightarrow~|p(x)-p_\infty|\leq\frac{C_\infty}{\log (e+|x|)}
\end{align*}
where $ p_\infty=\lim_{x\rightarrow\infty}p(x) $.
\end{itemize}
We will denote $ p(\cdot)\in LH(\mathbb{R}^n) $, if $ p(\cdot) $ is log-H\"{o}lder continuous locally and at infinity. 
\subsection{Variable Discrete Lebesgue Spaces}
For a sequence of exponents $ (p_n)_n $, we use the following notations
\begin{align*}
\bar{p}_n=\sup_n p_n\quad\&\quad\underline{p}_n=\inf_n p_n
\end{align*}
\begin{definition}
We define variable discrete Lebesgue space for a sequence of exponents $ (p_n)_n $ as follows 
\begin{align*}
&\ds:=\{(a_n)_n |~\|(a_n)_n\|_{\ds}<\infty\}
\\
&\quad s.t.\quad \|(a_n)_n\|_{\ds}=\inf\{\lambda>0:~\sum_{n=-\infty}^{\infty}|\frac{a_n}{\lambda}|^{p_n} \leq 1\}
\end{align*}
\end{definition}

\begin{definition}\label{def2}
The discrete Hilbert transform is an operator $\mathcal{H}$ that acts on a sequence $b = (b_n)_n$ by defining the sequence $\mathcal{H}b = (\mathcal{H}b_n)_n$, where for each $n \in \mathbb{Z}$,
\begin{align}\label{Hil}
\mathcal{H}b_n = \sum_{m \neq n} \frac{b_m}{n - m}.
\end{align}
\end{definition}

The operator defined in \eqref{Hil} was first introduced by D. Hilbert. It was shown by M. Riesz \cite{Riesz 30, Riesz 31} and E. C. Titchmarsh \cite{Titchmarsh 32, Titchmarsh 35} that for a constant $p \in (1, \infty)$, this operator is bounded on $\ell^p$. The purpose of this paper is to investigate the boundedness of the discrete Hilbert transform on $\ell^{p_n}(\mathbb{Z})$.

Assume that $\bar{p}_n < \infty$. For $b \in \ds$ and every $n \in \mathbb{Z}$, we have:
\begin{align*}
|\mathcal{H}b_n| &\leq \sum_{m \neq n} \frac{|b_m|}{|n - m|}
= \sum_{j = 1}^{\infty} \sum_{2^{j - 1} \leq |n - m| < 2^j} \frac{|b_m|}{|n - m|}
\\
&\leq \sum_{j = 1}^{\infty} \frac{1}{2^{j - 1}} \left(\sum_{|n - m| < 2^j} |b_m|^{\bar{p}_n} \right)^{\frac{1}{\bar{p}_n}} \left(2^{j + 1} - 1\right)^{1 - \frac{1}{\bar{p}_n}}
\\
&\leq \sum_{j = 1}^{\infty} \frac{\left(2^{j + 1}\right)^{1 - \frac{1}{\bar{p}_n}}}{2^{j - 1}} \left(\sum_{m} |b_m|^{\bar{p}_n} \right)^{\frac{1}{\bar{p}_n}}
\\
&= 4 \sum_{j = 1}^{\infty} \left(2^{-\frac{1}{\bar{p}_n}}\right)^{j + 1} \|b\|_{\ell^{\bar{p}_n}(\mathbb{Z})}
\\
&= \frac{4 \cdot 2^{-\frac{2}{\bar{p}_n}}}{1 - 2^{-\frac{1}{\bar{p}_n}}} \|b\|_{\ell^{\bar{p}_n}(\mathbb{Z})}.
\end{align*}
Since for every $n$, $p_n \leq \bar{p}_n$, it follows that $\|\cdot\|_{\ell^{\bar{p}_n}(\mathbb{Z})} \leq c \|\cdot\|_{\ds}$. Therefore,
\begin{align*}
|\mathcal{H}b_n| \leq c' \|b\|_{\ds} < \infty.
\end{align*}
Hence, Definition \ref{def2} is well-defined.

We now state the main result in the following theorem.

\begin{theorem}\label{theo:boundedness}
Let $\bar{p}_n < \infty$. Then, the discrete Hilbert transform $\mathcal{H}$ is bounded on $\ds$, i.e., 
\begin{align*}
\exists c(p) \quad \text{such that} \quad \|\mathcal{H}\|_{\ds \rightarrow \ds} \leq c(p).
\end{align*}
\end{theorem}

\section{Auxiliary lemmas and the proof of main result}

Before proving the theorem, we need to establish some lemmas. To demonstrate our main result, we introduce the following quantities. Consider the sequences $(b_k)_{k \in \mathbb{Z}}$ and $(p_k)_{k \in \mathbb{Z}}$. For each \( k \in \mathbb{Z} \), we define the functions \( f(x) \) and \( p(x) \) as follows:
\begin{align} \label{Yar}
f(x) := \left\{
\begin{array}{ll}
2\pi b_k &\text{if } |x - k| \leq \frac{1}{4}, \\
0 &\text{otherwise},
\end{array}
\right.
\quad \text{and} \quad
p(x) := \left\{
\begin{array}{ll}
p_k &\text{if } |x - k| \leq \frac{1}{2}, \\
0 &\text{otherwise}.
\end{array}
\right.
\end{align}

\begin{lemma}\label{lem:finLp()}
Let $\bar{p}_n < \infty$ and let $f$ and $p(\cdot)$ be as in \eqref{Yar}. Then, $f \in L^{p(\cdot)}(\mathbb{R})$.
\end{lemma}

\begin{proof}
Assume $\|b\|_{\ds} \leq 1$, $r \in (0, \infty)$, and $x_0 \in \mathbb{R}$. Then:
\begin{align*}
\forall r > 0: \exists n_0 \in \mathbb{N} \quad \text{such that} \quad r \in (n_0, n_0 + 1],
\\
\forall x_0 \in \mathbb{R}: \exists k_0 \in \mathbb{Z} \quad \text{such that} \quad |x_0 - k_0| \leq \frac{1}{2}.
\end{align*}
Therefore,
\begin{align*}
\varrho_{L^{p(\cdot)}(B(x_0,r))}(f)
&= \int_{x_0 - r}^{x_0 + r} |f(x)|^{p(x)} \, dx
\\
&\leq \int_{k_0 - (n_0 + 1) - \frac{1}{2}}^{k_0 + (n_0 + 1) + \frac{1}{2}} |f(x)|^{p(x)} \, dx
\\
&\leq \int_{k_0 - (n_0 + 1) - \frac{1}{4} - \frac{1}{4}}^{k_0 + (n_0 + 1) + \frac{1}{4} + \frac{1}{4}} |f(x)|^{p(x)} \, dx
\\
&= \sum_{|k - k_0| \leq n_0 + 1 + \frac{1}{4}} \int_{k - \frac{1}{4}}^{k + \frac{1}{4}} |f(x)|^{p(x)} \, dx
\\
&= \frac{1}{2} \sum_{|k - k_0| \leq n_0 + 1 + \frac{1}{4}} |2\pi b_k|^{p_k}
\\
&\leq \frac{(2\pi)^{\bar{p}_n}}{2} \sum_{k \in \mathbb{Z}} |b_k|^{p_k}
\\
&\leq \frac{(2\pi)^{\bar{p}_n}}{2} \|b\|_{\ds}.
\end{align*}
As $r \to \infty$, it follows that $\varrho_{L^{p(\cdot)}(\mathbb{R})}(f) < \infty$. Since $\esssup p(\cdot) = \bar{p}_n < \infty$, we conclude that $f \in L^{p(\cdot)}(\mathbb{R})$.
\end{proof}

\begin{theorem}\label{theo:HfinLp()}
\cite[Theorem 5.39]{Cruz} Given $p(\cdot) \in \mathcal{P}$ such that $1 < p_- \leq p^+ < \infty$, if the maximal operator is bounded on $L^{p(\cdot)}(\mathbb{R}^n)$, i.e.,
\begin{align*}
\frac{1}{p(\cdot)} \in LH \quad \text{and} \quad p_- > 1,
\end{align*}
then the Hilbert transform $\mathcal{H}$ is bounded on $L^{p(\cdot)}(\mathbb{R})$.
\end{theorem}

According to Theorem \ref{theo:HfinLp()}, for $ f\in L^{p(\cdot)}(\mathbb{R}^n) $, we have $ \mathcal{H}f\in L^{p(\cdot)}(\mathbb{R}) $. Let us define the following function which would help us to prove the Theorem \ref{theo:boundedness}. We consider the sequence $ (\mathcal{H}{b}_m)_{m\in\mathbb{Z}} $, then for $ n\in\mathbb{Z} $ we may define,   
\begin{align}\label{Yar-1}
F(x):=\left\{
\begin{array}{ll}
\mathcal{H}{b}_n &: ~|x-n|\leq\frac{1}{2}
\\
0 &: ~otherwise
\end{array}
\right.
\quad\&\quad
G(x):=\left\{
\begin{array}{ll}
\mathcal{H}f(x)-F(x) &: ~|x-n|\leq\frac{1}{2}
\\
0 &: ~otherwise
\end{array}
\right.
\end{align}

We now proceed to show that \( G \in L^{p(\cdot)}(\mathbb{R}) \).

\begin{lemma}\label{lem:GinLp()}
Let \( \bar{p}_n<\infty \) and \( G \) be as defined in \eqref{Yar-1}. Then, \( G \in L^{p(\cdot)}(\mathbb{R}) \).
\end{lemma}

\begin{proof}
Let \( \|b\|_{\ds}\leq 1 \), \( n\in\mathbb{Z} \), and \( |x-n|\leq \frac{1}{2} \) with \( x\neq n\pm\frac{1}{4} \). Then, we have
\begin{align*}
G(x) &= \frac{1}{\pi}\int_{\mathbb{R}}\frac{f(t)}{x-t}\,dt - \mathcal{H}{b}_n \\
&= \sum_{m\neq n} 2 b_m \int_{m-\frac{1}{4}}^{m+\frac{1}{4}} \left(\frac{1}{x-t} - \frac{1}{n-m}\right) \, dt + 2b_n \int_{n-\frac{1}{4}}^{n+\frac{1}{4}} \frac{1}{x-t} \, dt \\
&= G_1(x) + G_2(x).
\end{align*}

\begin{itemize}
\item[Case 1.] We first show that \( G_1 \in L^{p(\cdot)}(\mathbb{R}) \). For \( m\in\mathbb{Z} \) and \( m\neq n \), we have
\end{itemize}
\begin{align*}
\left\{
\begin{array}{ll}
|x-n|\leq\frac{1}{2}, \\
|m-t|\leq\frac{1}{4},
\end{array}
\right.
\Rightarrow |x-t|\geq |n-m|-\frac{3}{4} \Rightarrow \left|\frac{1}{x-t} - \frac{1}{n-m}\right| \leq \frac{3}{|n-m|^2}.
\end{align*}
Thus,
\begin{align*}
|G_1(x)| \leq \sum_{m\neq n} 2|b_m| \int_{m-\frac{1}{4}}^{m+\frac{1}{4}} \left|\frac{1}{x-t} - \frac{1}{n-m}\right| dt \leq \sum_{m\neq n} \frac{3|b_m|}{|n-m|^2}.
\end{align*}
Then, as in Lemma \ref{lem:finLp()}, we have
\begin{align*}
\varrho_{L^{p(\cdot)}(B(x_0,r))}(G_1) &= \int_{x_0-r}^{x_0+r} |G_1(x)|^{p(x)} \, dx \\
&\leq \int_{k_0-(n_0+1)-\frac{1}{2}}^{k_0+(n_0+1)+\frac{1}{2}} |G_1(x)|^{p(x)} \, dx \\
&= \sum_{|n-k_0|\leq n_0+1} \int_{n-\frac{1}{2}}^{n+\frac{1}{2}} |G_1(x)|^{p(x)} \, dx \\
&\leq 3^{\bar{p}_n} \sum_{|n-k_0|\leq n_0+1} \left|\sum_{m\neq n} \frac{|b_m|}{|n-m|^2} \right|^{p_n} \\
&\leq 3^{\bar{p}_n} \sum_{|n-k_0|\leq n_0+1} \left( \sum_{m\neq n} \frac{|b_m|^{\bar{p}_m}}{|n-m|^2} \right)^{\frac{1}{\bar{p}_m}} \left( \sum_{m\neq n} \frac{1}{|n-m|^2} \right)^{\frac{1}{\bar{p}'_m}} 
\\
&\leq 3^{\bar{p}_n} 4^{\frac{\bar{p}_n}{\bar{p}'_m}} \sum_{|n-k_0|\leq n_0+1} \left( \sum_{m\neq n} \frac{|b_m|^{\bar{p}_m}}{|n-m|^2} \right)^{\frac{p_n}{\bar{p}_m}} \\
&\leq 3^{\bar{p}_n} 4^{\frac{\bar{p}_n}{\bar{p}'_m}} \sum_{|n-k_0|\leq n_0+1} \left( \sum_{m\neq n} \frac{|b_m|^{p_m}}{|n-m|^2} \right)^{\frac{\underline{p}_n}{\bar{p}_m}} 
\end{align*}
\begin{align*}
&\leq 3^{\bar{p}_n} 2^{2\frac{\bar{p}_n}{\bar{p}'_m}} \max\{1,2^{1-\frac{\underline{p}_n}{\bar{p}_m}}\} \left( \sum_{|n-k_0|\leq n_0+1} \sum_{m\neq n} \frac{|b_m|^{p_m}}{|n-m|^2} \right)^{\frac{\underline{p}_n}{\bar{p}_m}} \\
&= 3^{\bar{p}_n} 2^{1+\frac{\underline{p}_n}{\bar{p}_m}} \left( \sum_{m} \sum_{|n-k_0|\leq n_0+1,~n\neq m} \frac{|b_m|^{p_m}}{|n-m|^2} \right)^{\frac{\underline{p}_n}{\bar{p}_m}} \\
&\leq 12^{\bar{p}_n} \left( \sum_{|m-k_0|\leq 2(n_0+1)} |b_m|^{p_m} 4 \right. \\
&\left.+ \sum_{i=1}^{\infty} \sum_{~2^i(n_0+1)<|m-k_0|\leq 2^{i+1}(n_0+1)~} \sum_{|n-k_0|\leq n_0+1,~n\neq m} \frac{|b_m|^{p_m}}{|n-m|^2} \right)^{\frac{\underline{p}_n}{\bar{p}_m}}.
\end{align*}
Since
\begin{align*}
2^i(n_0+1) \leq |m-k_0| &\leq |n-m| + |n-k_0| \leq |n-m| + n_0+1,
\end{align*}
then
\begin{align*}
\left\{
\begin{array}{ll}
(n_0+1)(2^i-1) \leq |n-m|, \\
(n_0+1)^{\frac{1}{2}}(2^{i-\frac{3}{2}}) \leq (n_0+1)(2^i-1),
\end{array}
\right.
\Rightarrow (n_0+1)^{\frac{1}{2}}(2^{i-\frac{3}{2}}) \leq |n-m|.
\end{align*}
Therefore,
\begin{align*}
\varrho_{L^{p(\cdot)}(B(x_0,r))}(G_1) &\leq 12^{\bar{p}_n} \left( 4\|b\|_{\ds} + \sum_{i=1}^{\infty} \sum_{m} \sum_{|n-k_0|\leq n_0+1,~n\neq m} \frac{|b_m|^{p_m}}{(n_0+1)2^{2i-3}} \right)^{\frac{\underline{p}_n}{\bar{p}_m}} \\
&\leq 12^{\bar{p}_n} \left( 4\|b\|_{\ds} + \sum_{m} |b_m|^{p_m} \sum_{i=1}^{\infty} 2(n_0+1) \frac{1}{(n_0+1)2^{2i-3}} \right)^{\frac{\underline{p}_n}{\bar{p}_m}} \\
&\leq 12^{\bar{p}_n} \left( 4\|b\|_{\ds} + \sum_{m} |b_m|^{p_m} \sum_{i=1}^{\infty} \frac{4}{2^{2i-3}} \right)^{\frac{\underline{p}_n}{\bar{p}_m}} \\
&= 12^{\bar{p}_n} \left( 4\|b\|_{\ds} + \sum_{m} |b_m|^{p_m} \right)^{\frac{\underline{p}_n}{\bar{p}_m}} \\
&\leq 12^{\bar{p}_n} 4^{\frac{\underline{p}_n}{\bar{p}_m}} (4\|b\|_{\ds} + 1)^{\frac{\underline{p}_n}{\bar{p}_m}} \leq 12^{\bar{p}_n} 4^{1+\frac{\underline{p}_n}{\bar{p}_m}}.
\end{align*}
\item[Case 2.] Now, we show that \( G_2 \in L^{p(\cdot)}(\mathbb{R}) \). Let \( x=n+t \) and \( -\frac{1}{2}\leq t\leq \frac{1}{2} \). We have
\begin{align*}
G_2(x) &= 2b_n\int_{-\frac{1}{4}}^{\frac{1}{4}} \frac{1}{t} dt \\
&= 2b_n \log \left|\frac{n+\frac{1}{4}}{n-\frac{1}{4}}\right| \leq 2|b_n| \log 2.
\end{align*}
Thus,
\begin{align*}
\varrho_{L^{p(\cdot)}(B(x_0,r))}(G_2) &= \int_{x_0-r}^{x_0+r} |G_2(x)|^{p(x)} dx = \sum_{|n-k_0|\leq n_0+1} \int_{n-\frac{1}{2}}^{n+\frac{1}{2}} |G_2(x)|^{p(x)} dx \\
&\leq 2^{\bar{p}_n} (\log 2)^{\bar{p}_n} \sum_{|n-k_0|\leq n_0+1} |b_n|^{p_n}.
\end{align*}
Finally, we have \( G(x) = G_1(x) + G_2(x) \), and \( G \in L^{p(\cdot)}(\mathbb{R}) \).
\end{proof}

\begin{lemma}\label{FinLp()}
Let $ 1<p_-\leq p^+<\infty $ and $ 1\slash p(\cdot)\in LH $ and $ F $ be as in \eqref{Yar-1}, then $ F\in L^{p(\cdot)}(\mathbb{R}) $.
\end{lemma}
\begin{proof}
Let $ \|b\|_{\ds}\leq 1 $. For $ k_0\in\mathbb{Z} $ and $ n_0\in\mathbb{N}_0 $ we have we have
\begin{align*}
\varrho_{L^{p(\cdot)}(B(k_0,(n_0+1)+\frac{1}{2}))}(F)&=\int_{k_0-(n_0+1)-\frac{1}{2}}^{k_0+(n_0+1)+\frac{1}{2}}|F(x)|^{p(x)}dx
\\
&=\sum_{|n-k_0|\leq n_0+1}\int_{n-\frac{1}{2}}^{n+\frac{1}{2}}|F(x)|^{p(x)}dx,
\end{align*}
according to definition of $ F $, Theorem \ref{theo:HfinLp()} and Lemma \ref{lem:GinLp()} we have
\begin{align*}
\sum_{|n-k_0|\leq n_0+1}|\mathcal{H}{b}_n|^{p_n}&=\varrho_{L^{p(\cdot)}(B(k_0,(n_0+1)+\frac{1}{2}))}(F)
\\
&\leq 2^{\bar{p}_n-1}\Big[\varrho_{L^{p(\cdot)}(B(k_0,(n_0+1)+\frac{1}{2}))}(G)+\varrho_{L^{p(\cdot)}(B(k_0,(n_0+1)+\frac{1}{2}))}(\mathcal{H}f)\Big]
\\
&\leq 2^{\bar{p}_n-1}\Big[\varrho_{L^{p(\cdot)}(\mathbb{R})}(G)+\varrho_{L^{p(\cdot)}(\mathbb{R})}(\mathcal{H}f)\Big]
\\
&<\infty.
\end{align*}
Passing $ n_0\rightarrow\infty $, we have the desired result.
\end{proof}
Now we refer to the proof of the Theorem \ref{theo:boundedness}.
\begin{proof}[{\bf Proof of Theorem \ref{theo:boundedness}}]
Let $ \|b\|_{\ell^{p_n}(\mathbb{Z})}\leq 1 $. According to Lemma \ref{FinLp()} we have
\begin{align*}
\exists M\geq 1\quad:\quad\varrho_{\ell^{p_n}(\mathbb{Z})}(\mathcal{H}{b})=\varrho_{L^{p(\cdot)}(\mathbb{R})}(F)\leq M,
\end{align*}
therefore $ \|\mathcal{H}{b}\|_{\ell^{p_n}(\mathbb{Z})}\leq M $, which means
\begin{align*}
\|\mathcal{H}\|_{\ds\rightarrow\ds}&=\sup_{\|b\|_{\ds}\leq 1}\|\mathcal{H}b\|_{\ds}\leq\quad M
\end{align*}
which is the desired result.
\end{proof}

\section{The Discrete Mikhlin Multiplier}

Let $ \mathcal{F}f $ stands for the Fourier transform of $ f $. We investigate the Mikhlin multiplier property in the variable discrete Lebesgue space $ \ell^{p_n}(\mathbb{Z}) $.

\begin{definition}\label{def5}
The discrete Mikhlin multiplier is an operator $ \mathcal{T}_m $ such that for a sequence $ b=(b_n)_n $ and $ m\in L^\infty $, we define $ \mathcal{T}_m b=(T_mb_n)_n $ where for every $ n\in\mathbb{Z} $  we have
\begin{align*}
T_mb_n=\mathcal{F}^{-1}m*b_n.
\end{align*}
\end{definition}

Let $ b\in\ds, m\in L^\infty, n\in\mathbb{Z} $ and $ y\in\mathbb{R} $ then
\begin{align*}
|T_mb_n(y)|&=|(\mathcal{F}^{-1}m*b_n)(y)|
\\
&=|b_n\int_{\mathbb{R}}\mathcal{F}^{-1}m(x)dx|
\\
&=|b_n|~|\int_{\mathbb{R}}\int_{\mathbb{R}}m(\xi)e^{i2\pi \xi x}d\xi dx|
\\
&=|b_n|~|\int_{\mathbb{R}}m(\xi)\int_{\mathbb{R}}e^{i2\pi \xi x} dx~d\xi|
\\
&=|b_n|~|\int_{\mathbb{R}}m(\xi)~(\mathcal{F}^{-1} 1)(\xi)~d\xi|
\\
&=|b_n|~|\int_{\mathbb{R}}m(\xi)~\delta(\xi)~d\xi|
\\
&\leq |b_n|~\|m\|_{L^\infty}~\int_{\mathbb{R}}\delta(\xi)~d\xi
\end{align*}
since for Dirac delta function we have $ \|\delta\|_{L^1}=1 $, then
\begin{align*}
|T_mb_n(y)|\leq |b_n|~\|m\|_{L^\infty}<\infty.
\end{align*}
Therefore, the Definition \ref{def5} is well defined. 

Now we aim to demonstrate that the discrete Mikhlin multiplier operator \( \mathcal{T}_m \) is bounded on \( \ell^{p_n}(\mathbb{Z}) \).

\begin{remark}\label{rem:modular}
Given the definition of the modular and equation \eqref{Yar}, we have:
\begin{align*}
\varrho_{L^{p(\cdot)}}(f) &= \int_{\mathbb{R}} |f(x)|^{p(x)} \, dx = \sum_{n\in\mathbb{Z}} \int_{B(n,\frac{1}{4})} |f(x)|^{p(x)} \, dx = \frac{1}{2}\sum_{n\in\mathbb{Z}} |2\pi b_n|^{p_n}.
\end{align*}
Therefore, there exist constants \( m_1, M_1 \in [1, \infty) \) such that:
\begin{align*}
m_1 ~ \varrho_{\ell^{p_n}(\mathbb{Z})}(b) \leq \varrho_{L^{p(\cdot)}(\mathbb{R})}(f) \leq M_1 ~ \varrho_{\ell^{p_n}(\mathbb{Z})}(b).
\end{align*}
\end{remark}

\begin{theorem}\label{theo1}
Let \( f \) and \( p(\cdot) \) be as defined in equation \eqref{Yar}, with \( 1 < p_- \leq p_+ \leq 2 \), and suppose \( 1 \in L^{r(\cdot)}(D) \), where 
\begin{align*}
D = \{x : p_- < p(x)\}, \qquad \text{and} \qquad \frac{1}{p_-} = \frac{1}{p(\cdot)} + \frac{1}{r(\cdot)}.
\end{align*}
Let \( m : \mathbb{R}^d \rightarrow \mathbb{C} \) be such that
\begin{align*}
\mathcal{F}^{-1} m \in L^q, \qquad \text{and} \qquad \frac{1}{p_-} + \frac{1}{q} = 1 + \frac{1}{2}.
\end{align*}
Define the operator \( T_m \) by
\begin{align*}
T_m f = \mathcal{F}^{-1}(m~ \mathcal{F} f),
\end{align*}
which is a linear \( L^2 \)-bounded and (1,1)-weak type operator. Then, the operator \( \mathcal{T}_m \) is bounded on \( \ell^{p_n}(\mathbb{Z}) \) for \( 1 < p_n < 2 \). That is, there exists a constant \( M \) such that
\begin{align*}
\|\mathcal{T}_m\|_{\ell^{p_n}(\mathbb{Z}) \rightarrow \ell^{p_n}(\mathbb{Z})} < M.
\end{align*}
\end{theorem}

\begin{proof}
Let \( \|b\|_{\ell^{p_n}(\mathbb{Z})} \leq 1 \). Since \( \bar{p}_n < \infty \), by the linearity of \( \mathcal{T}_m \) and by \cite[Proposition 2.12]{Cruz}, it suffices to show that \( \varrho_{\ell^{p_n}(\mathbb{Z})}(\mathcal{T}_m b) < \infty \).

By the Marcinkiewicz interpolation theorem (note that \( 1 < p_- \)) and Young's inequality \cite[Proposition 5.2]{Cruz}, combined with the embedding theorem \cite[Theorem 2.45]{Cruz}, we have:
\begin{align*}
\varrho_{\ell^{p_n}(\mathbb{Z})}(\mathcal{T}_m b) &\leq m^{-1}_1 \varrho_{L^{p(\cdot)}(\mathbb{R})}(T_m f) \\
&= m^{-1}_1 \int_{\mathbb{R}^n} |T_m f(x)|^{p(x)} \, dx \\
&= m^{-1}_1 \left( \int_{\{x : |T_m f(x)| \leq 1\}} |T_m f(x)|^{p(x)} \, dx + \int_{\{x : |T_m f(x)| > 1\}} |T_m f(x)|^{p(x)} \, dx \right) \\
&\leq m^{-1}_1 \left( \int_{\{x : |T_m f(x)| \leq 1\}} |T_m f(x)|^{p_-} \, dx + \int_{\{x : |T_m f(x)| > 1\}} |T_m f(x)|^{2} \, dx \right) \\
&\leq m^{-1}_1 \left( \|T_m f\|_{p_-}^{p_-} + \|T_m f\|_2^2 \right) \\
&\leq m^{-1}_1 \left( c_1 \|f\|_{p_-}^{p_-} + \|\mathcal{F}^{-1} m * f\|_2^2 \right) \\
&\leq m^{-1}_1 \left( c_1 \|f\|_{p_-}^{p_-} + \|\mathcal{F}^{-1} m\|_q^2 \|f\|_{p_-}^2 \right) \\
&\leq m^{-1}_1 \left( c_1 K \|f\|_{p(\cdot)}^{p_-} + \|\mathcal{F}^{-1} m\|_q^2 K \|f\|_{p(\cdot)}^2 \right) \\
&\leq m^{-1}_1 \left( c_1 K M_1' \|b\|_{\ell^{p_n}(\mathbb{Z})}^{p_-} + \|\mathcal{F}^{-1} m\|_q^2 K M_1' \|b\|_{\ell^{p_n}(\mathbb{Z})}^2 \right) \\
&\leq m^{-1}_1 \left( c_1 K M_1' + K M_1' \|\mathcal{F}^{-1} m\|_q^2 \right) \\
&< \infty,
\end{align*}
which is the desired result.
\end{proof}

\begin{theorem}
Let \( f \) and \( p(\cdot) \) be as defined in the equation \eqref{Yar}, with \( p_+ < \infty \), and suppose \( 1 \in L^{r(\cdot)}(D) \), where 
\begin{align*}
D = \{x : p_- < p(x)\}, \qquad \text{and} \qquad \frac{1}{p_-} = \frac{1}{p(\cdot)} + \frac{1}{r(\cdot)}.
\end{align*}
Let \( m : \mathbb{R}^d \rightarrow \mathbb{C} \) be such that:
\begin{align*}
\mathcal{F}^{-1} m \in L^q, \qquad \text{and} \qquad \frac{1}{p_-} + \frac{1}{q} = 1 + \frac{1}{2},
\end{align*}
and for every \( \xi \neq 0 \), 
\begin{align}\label{MMC}
|\xi|^{|\gamma|} |\partial^{\gamma} m(\xi)| < B,
\end{align}
where \( \gamma \) is a multi-index with \( |\gamma| \leq d + 2 \). Then there exists a constant \( c > 0 \) such that:
\begin{align*}
\|\mathcal{T}_m\|_{\ell^{p_n}(\mathbb{Z}) \rightarrow \ell^{p_n}(\mathbb{Z})} < c.
\end{align*}
\end{theorem}

\begin{proof}
Let $ \|b\|_{\ell^{p_n}(\mathbb{Z})}\leq 1 $ and $ T_mf:=\mathcal{F}^{-1}(m~\mathcal{F}f) $.
For $ 2\leq p(\cdot) $ we will have $ 1<p'(\cdot)<2 $. On the other side, $ m\in L^\infty $, then by Plancherel theorem $ T_m $ is a linear $ L^2 $-bounded operator. Therefore, $ T_m $ is a Calderon-Zygmund operator such that by \eqref{MMC}, it satisfies in (gradient) Hormander condition, then it is (1,1)-weak type operator too. Therefore by Theorem \ref{theo1} there exists $ M_2>0 $, such that
\begin{align}\label{eq1:theo1}
\underset{\|g\|_{p'(\cdot)}\leq 1}{\sup}~ \|T_{m}g\|_{p'(\cdot)}\leq M_2\qquad s.t.\qquad 1<p'(\cdot)<2.
\end{align}
Note that adjoint kernel is $ m^*(x)=m(-x) $, then by associated norm \cite[Theorem 2.34]{Cruz} and Holder inequality \cite[Theorem 2.26]{Cruz}, for $ 2\leq p(\cdot) $ we have
\begin{align*}
\|T_mf\|_{p(\cdot)}&\leq k_{p(\cdot)}^{-1}~\underset{\|g\|_{p'(\cdot)}\leq 1}{\sup}\int |(T_mf)(x)g(x)|
\\
&\leq k_{p(\cdot)}^{-1}~\underset{\|g\|_{p'(\cdot)}\leq 1}{\sup}\int |f(x)(T_{m^*}g)(x)|
\\
&\leq k_{p(\cdot)}^{-1}~\underset{\|g\|_{p'(\cdot)}\leq 1}{\sup}K_{p(\cdot)}~\|f\|_{p(\cdot)} \|T_{m^*}g\|_{p'(\cdot)}
\\
&=k_{p(\cdot)}^{-1}~K_{p(\cdot)}\|f\|_{p(\cdot)}\underset{\|g\|_{p'(\cdot)}\leq 1}{\sup}~ \|T_{m}g\|_{p'(\cdot)}
\end{align*}
therefore, by equation \eqref{eq1:theo1}, we have
\begin{align}\label{eq:bounded}
\|T_mf\|_{p(\cdot)}\leq \Big(k_{p(\cdot)}^{-1}~K_{p(\cdot)}~M_2\Big)\|f\|_{p(\cdot)}\qquad s.t.\qquad 2\leq p(\cdot).
\end{align}
Since $ \bar{p}_n<\infty $ and $ \mathcal{T}_m $ is a linear operator, then by \cite[Proposition 2.12]{Cruz}, it suffices to show $ \varrho_{\ell^{p_n}(\mathbb{Z})}(\mathcal{T}_mb)<\infty $. Now by considering \eqref{eq1:theo1} and \eqref{eq:bounded} and $ p_+<\infty $, as in the Remark \ref{rem:modular} we have
\begin{align*}
\varrho_{\ell^{p_n}(\mathbb{Z})}(\mathcal{T}_mb)& \leq m^{-1}_1\varrho_{p(\cdot)}(T_mf)
\\
&=m^{-1}_1\int_{\mathbb{R}^n}|T_mf(x)|^{p(x)}
\\
&=m^{-1}_1\Big(\int_{\{x:1<p(x)<2\}}|T_mf(x)|^{p(x)}+\int_{\{x:2\leq p(x)\}}|T_mf(x)|^{p(x)}\Big)
\\
&<\infty,
\end{align*}
which completes the proof.
\end{proof}

\subsection*{Acknowledgements} {}

\bibliographystyle{amsplain}

\end{document}